\input amstex  
\magnification =\magstep 1
\documentstyle{amsppt}
\pageheight{9truein}
\pagewidth{6.5truein}
\NoRunningHeads
\baselineskip=16pt

\topmatter

\title The second maximal groups with respect to the sum of element orders
\endtitle

\author Marcel Herzog*, Patrizia Longobardi** and Mercede Maj**
\endauthor

\affil *School of Mathematical Sciences \\
       Tel-Aviv University \\
       Ramat-Aviv, Tel-Aviv, Israel
{}\\
               **Dipartimento di Matematica \\
       Universit\`a di Salerno\\
       via Giovanni Paolo II, 132, 84084 Fisciano (Salerno), Italy
\endaffil
\abstract Denote by $G$ a finite group and let $\psi(G)$ denote the sum of
element orders in $G$. In 2009, H.Amiri, S.M.Jafarian Amiri and I.M.Isaacs
proved that if $|G|=n$ and $G$ is non-cyclic, then $\psi(G)<\psi(C_n)$, where $C_n$ denotes the
cyclic group of order $n$.  In 2018 we proved that if $G$ is non-cyclic group
of order $n$, then $\psi(G)\leq \frac 7{11}\psi(C_n)$ and equality holds if
$n=4k$ with $(k,2)=1$ and $G=(C_2\times C_2)\times C_k$. In this paper we 
proved that equality holds if and only if $n$ and $G$ are as indicated above.
Moreover we proved the following generalization of this result: Theorem 4.
Let $q$ be a prime and let $G$ be a non-cyclic group of order $n$, 
with $q$ being the least prime divisor of $n$. Then $\psi(G)\leq 
\frac {((q^2-1)q+1)(q+1)}{q^5+1}\psi(C_n)$, with equality if and only if
$n=q^2k$ with $(k,q)=1$ and $G=(C_q\times C_q)\times C_k$. Notice that if $q=2$,
then  $\frac {((q^2-1)q+1)(q+1)}{q^5+1}=\frac 7{11}$.
\endabstract

\thanks This work was supported by the National Group for Algebraic and
Geometric Structures, and their Applications (GNSAGA - INDAM), Italy.
The first author  is grateful to the Department of
Mathematics of the University of Salerno
for its hospitality and
support, while this investigation was carried out.
\endthanks
\endtopmatter

\document
\heading I. Introduction \\
\endheading

In this paper all groups are finite. In [1], H. Amiri,
S.M. Jafarian Amiri and I.M. Isaacs introduced the following function on groups:
$$\psi(G)=\sum\{o(x)\ \mid\ x\in G\},$$
where $o(x)$ denotes the order of $x$. Thus $\psi(G)$
denotes the {\it sum}
of  element  orders  of  the  finite group $G$. Denoting the
cyclic group of order $n$ by $C_n$, they proved the following theorem:
\proclaim {Theorem AAI} If $G$ is a {\it non-cyclic} finite group of order $n$,
then
$$\psi(G)\ <\ \psi(C_n).$$
\endproclaim
Thus the  group $C_n$ is the unique group of order $n$ which attains the maximal
value of $\psi(G)$ among groups of order $n$. We shall call the groups $C_n$
"maximal" groups.

Starting from this results recently many authors have studied the function $\psi(G)$ and its relations with the structure of $G$ (see for example [2]-[17]).
 In the papers [8] and [16] M. Amiri and S.M. Jafarian Amiri, and, independently, R. Shen, G. Chen and C. Wu started the investigation of groups with the second largest value of the sum of element orders.

In [12], we determined the exact upper bound for $\psi(G)$ for non-cyclic
groups of order $n$. We proved the following theorem:
\proclaim{Theorem 1}
If $G$ is a non-cyclic group of order $n$, then
$$\psi(G)\ \leq\ \frac 7{11}\psi(C_n).$$
Moreover, the equality holds if $n=4k$ with $(k,2)=1$ and $G=(C_2\times
C_2)\times C_k$.
\endproclaim
Groups of order $n$ satisfying $\psi(G)\ =\ \frac 7{11}\psi(C_n)$
will be called "second maximal". Theorem 1
left open the following problem:
\proclaim{Problem} Determine all second maximal groups.
\endproclaim

In this paper we present a solution to this problem. We prove 
\proclaim {Theorem 2} The group $G$ of order $n$ is second maximal if and only
if $n=4k$ with $(k,2)=1$ and $G=(C_2\times C_2)\times C_k$.
\endproclaim
Theorem 1 and Theorem 2 imply the following complete result.
\proclaim {Theorem 3} If $G$ is a non-cyclic group of order $n$, then
$$\psi(G)\ \leq\ \frac 7{11}\psi(C_n).$$
Moreover, the equality holds if and only if $n=4k$ with $(k,2)=1$ 
and $G=(C_2\times
C_2)\times C_k$.
\endproclaim
 
Theorem 2 follows from the following more general result. First, we introduce
some notation and remarks. If $q$ is a prime, we 
shall call a group $G$ a "$q^*$-group" if $q$ is the 
{\bf least prime divisor} of
the order of $G$. We also define the following function of the real variable:
$$f(x)=\frac {((x^2-1)x+1)(x+1)}{x^5+1}.$$
This function is strictly decreasing for $x\geq 2$. Indeed, if $x\geq 2$, then
$$f(x)=\frac {((x^2-1)x+1)(x+1)}{x^5+1}=\frac {x^3-x+1}{x^4-x^3+x^2-x+1}$$
and
$$f'(x)=\frac {-x(x^3(x^2-4) +5x^2 +(x-2)(3x-1))}{(x^4-x^3+x^2-x+1)^2}<0.$$

Our general result is the following theorem.
\proclaim {Theorem 4} Let $q$ be a prime and let $G$ be a non-cyclic
$q^*$-group
of order $n$.
Then
$$\psi(G)\leq \frac {((q^2-1)q+1)(q+1)}{q^5+1}\psi(C_n)=f(q)\psi(C_n)$$
and the equality holds if and only if $n=q^2k$ with $(k,q!)=1$ and
$$G=(C_q\times C_q)\times C_k.$$
\endproclaim 

Notice that if $q=2$, then $f(q)=\frac {21}{33}=\frac 7{11}$. Moreover,
if $q$ is an odd prime, then $q>2$ and by our previous remark
$f(q)<f(2)=\frac 7{11}$. Hence it follows by Theorem 4 that if $G$
is a $q^*$-group for some odd prime $q$, then 
$$\psi(G)\leq f(q)\psi(C_n)<\frac 7{11}\psi(C_n)$$
and the group $G$ is not second maximal. Hence the second maximal groups
are $2^*$-groups $G$ satisfying the identity 
$\psi(G)=f(2)\psi(C_n)$ and by Theorem 4 these groups are the groups 
$G=(C_2\times C_2)\times C_k$ with $n=4k$ and $(k,2)=1$. Thus Theorem 2
follows from Theorem 4.

We shall call $q^*$-groups $G$ of order $n$ with $\psi(G)$
satisfying the identity $\psi(G)=f(q)\psi(C_n)$
"second maximal $q^*$-groups". If follows by Theorem 4 that
the
second maximal groups are exactly the second maximal $2^*$-groups, which
were
determined in Theorem 4. 

In addition to determining exactly the second maximal groups,
Theorem 4 also determines  exactly the second maximal $q^*$-groups 
for all primes $q$. This remark will be our last theorem.
\proclaim {Theorem 5} Let $q$ denote a prime. The $q^*$-group $G$ 
of order $n$ is a second maximal ${q^*}$-group if and only
if $n=q^2k$ with $(k,q!)=1$ and $G=(C_q\times C_q)\times C_k$.
\endproclaim

In order to conclude this paper, we need only to prove Theorem 4.
The proof of Theorem 4 is achieved in two steps. The first step 
is the  
following proposition:
\proclaim {Proposition 6} Let $q$ be a prime and let $G$ be a non-cyclic 
$q^*$-group
of order $n$. 
Then
$$\psi(G)\leq \frac {((q^2-1)q+1)(q+1)}{q^5+1}\psi(C_n).$$
\endproclaim
 
The second, and final step  in our proof of Theorem 4 is the following
proposition:
\proclaim {Proposition 7}  Let $q$ be a prime and let $G$ be a  
$q^*$-group of
order $n$. Then
$$\psi(G)= \frac {((q^2-1)q+1)(q+1)}{q^5+1}\psi(C_n)$$
if and only if $n=q^2k$ with $(k,q!)=1$ and
$$G=(C_q\times C_q)\times C_k.$$
\endproclaim
 
It is clear that Theorem 4 follows from Propositions 6 and 7. The proofs 
of these propositions will be presented in Sections II and III, respectively. 

\heading II. Proof of Proposition 6.
\endheading
  Before beginning with the proof, we shall list some preliminary results from 
the papers [1], [12] and [14].
\proclaim {Lemma 2.1}  Here $G$ denotes a finite group, $p,p_i$ denote primes
and $n,m,\alpha_i$ denote positive integer. The following statements hold.
\roster
\item"({\bf 1})" {\rm ([12], Lemma 2.9(1))} $\psi(C_{p^m})=\frac {p^{2m+1}+1}{p+1}=
\frac {p(p^{m})^2+1}{p+1}$.
\item"({\bf 2})" {\rm ([12], Lemma 2.2(3))} If $G=A\times B$, where $A,B$ 
are subgroups of $G$ satisfying
$(|A|,|B|)=1$, then
$\psi(G)=\psi(A)\psi(B)$.
\item"({\bf 3})" {\rm ([12], Lemma 2.9(2))} 
If $n=\prod _{i=1}^{i=m}p_i^{\alpha_i}$, 
where $p_i\neq p_j$ for
$i\neq j$, then
$\psi(C_n)=\prod _{i=1}^{i=m}\psi(C_{p_i^{\alpha_i}})$.
\item"({\bf 4})" {\rm ([12], Proof of Lemma 2.9(2))} 
$\psi(C_n)\geq \frac q{p+1}n^2$.
\item"({\bf 5})" {\rm ([1], Corollary B)} If $P$ is a cyclic normal Sylow 
subgroup of $G$ then 
$\psi(G)\leq \psi(P)\psi(G/P)$, with equality if and only if $P$ is 
central in $G$.
\item"({\bf 6})" {\rm ([12], Lemma 2.2(5))} Let $G=P\rtimes F$, where $P$ is 
a cyclic $p$-group, $|F|>1$, 
$(p,|F|)=1$ and $Z=C_F(P)$.
Then
$$ \psi(G)<\psi(P)\psi(F)\biggl(\frac {\psi(Z)}{\psi(F)}+\frac {|P|}{\psi(P)}
\biggr).$$ 
\item"({\bf 7})" {\rm ([14], Lemma 2.7)}  Let $G=P\rtimes F$, where $P$ is a cyclic Sylow 
$p$-subgroup of $G$ and $F$ 
is a cyclic $p$-complement. If $\psi(G)$ has the second largest value 
for groups of 
order $|G|$, then $C_F(P)$ is a maximal subgroup of $F$. 
\endroster
\endproclaim

We now begin with the proof of Proposition 6, which is restated below.
\proclaim {Proposition 6} Let $q$ be a prime and let $G$ be a non-cyclic $q^*$-group
of order $n$.
Then
$$\psi(G)\leq \frac {((q^2-1)q+1)(q+1)}{q^5+1}\psi(C_n).$$
\endproclaim
\demo {Proof}  
From now on, we shall assume that
$G$ is a non-cyclic $q^*$-group of order $n$ 
satisfying the inequality
$$\psi(G) > \frac {((q^2-1)q+1)(q+1)}{q^5+1}\psi(C_n),$$
and our aim is to reach a contradiction. Let $p$ denote the largest prime
divisor of $n$. Our proof is by induction on the size of $p$. Notice that by
our assumptions $n\geq q^2$.

If $q=2$, then
$\psi(G) > f(2)\psi(C_n)=
\frac 7{11} \psi(C_n)$
and since $G$ is non-cyclic, this contradicts Theorem 1.
Therefore we may assume that $q>2$.

By Lemma 2.1(4)
$\psi(C_n)\geq \frac q{p+1}n^2$, which implies that
$$\psi(G)> \frac {((q^2-1)q+1)(q+1)q}{(q^5+1)(p+1)}n^2.$$
Consequently, there exists $x\in G$ such that $o(x)>\frac
{((q^2-1)q+1)(q+1)q}{(q^5+1)(p+1)}n,$
which implies that
$$[G:\langle x\rangle]<\frac {(q^5+1)(p+1)}{((q^2-1)q+1)(q+1)q}.\tag{$*$}$$

Suppose, first, that $p=q$. Then  $n$ is a power of $q$ and by $(*)$
$$[G:\langle x\rangle]<\frac {q^5+1}{((q^2-1)q+1)q}.$$
Since
$((q^2-1)q+1)q>(q^2+1)q>\frac {q^5+1}{q^2}$, it follows that
$[G:\langle x\rangle]<q^2$. Thus
$[G:\langle x\rangle]=q$, which implies that
$$\psi(G)\leq \psi(C_{n/q})+(\frac {q-1}qn)\frac nq=\frac {q(n/q)^2+1}{q+1}+\frac
{q-1}{q^2}n^2=
\frac {q^2+q-1}{q^2(q+1)}n^2+\frac 1{q+1}.$$

We claim that
$\psi(G) \leq \frac {((q^2-1)q+1)(q+1)}{q^5+1}\psi(C_n)$,
yielding a contradiction to our assumption. It suffices to prove that
$$\frac {q^2+q-1}{q^2(q+1)}n^2+\frac 1{q+1}\leq \frac
{((q^2-1)q+1)(q+1)}{q^5+1}\cdot
\frac {qn^2+1}{q+1},$$
or multiplying by $(q+1)(q^5+1)q^2$, that
$$(q^2+q-1)(q^5+1)n^2+q^2(q^5+1)\leq (q^3-q+1)(q+1)q^3n^2+(q^3-q+1)(q+1)q^2.$$
Notice that $(q^2+q-1)(q^5+1)=q^7+q^6-q^5+q^2+q-1$ and
$(q^3-q+1)(q+1)q^3=q^7+q^6-q^5+q^3$.
Hence it suffices to prove that
$$(q^7+q^2)-(q^6+q^5-q^4+q^2)\leq (q^3-q^2-q+1)n^2.$$
But
$(q^7+q^2)-(q^6+q^5-q^4+q^2)=q^4(q^3-q^2-q+1)$,
so it suffices to prove that $q^4\leq n^2$, which is true since $n\geq q^2$. This
proves our claim and completes the proof in the case when $p=q$.

So suppose that $p>q$. Since $q>2$, it follows that $p\geq q+2$.
We claim that
$$[G:\langle x\rangle]<\frac {(q^5+1)(p+1)}{((q^2-1)q+1)(q+1)q}<p.$$
We need to prove that
$((q^3-q+1)(q^2+q)-(q^5+1))p>q^5+1$
or $(q^4-q^3+q-1)p>q^5+1$. But $p\geq q+2$ and $q\geq 3$, so
$$(q^4-q^3+q-1)p\geq (q^4-q^3+q-1)(q+2)>(q^4-q^3)(q+2)=q^5+q^4-2q^3>q^5+1,$$
as required. So our claim is true.

Hence $[G:\langle x\rangle]<p$ and $\langle x\rangle$ contains a cyclic 
Sylow $p$-subgroup $P$ of $G$. Since
$\langle x\rangle\leq N_G(P)$ and $[G:\langle x\rangle]<p$, it follows that
$P$ is a cyclic normal subgroup of $G$. Now our assumptions and Lemma 2.1(5) 
imply that
$$\psi(P)\psi(G/P)\geq \psi(G)>\frac
{((q^2-1)q+1)(q+1)}{q^5+1}\psi(C_{|P|})\psi(C_{|G/P|}).$$
Since $P\cong C_{|P|}$, cancellation yields
$$\psi(G/P)>\frac {((q^2-1)q+1)(q+1)}{q^5+1}\psi(C_{|G/P|}).$$
Since $G/P$ is a $q^*$-group and the maximal prime dividing 
$|G/P|$ is smaller than $p$, our inductive
hypothesis
implies that $G/P$ is a cyclic group and $G=P\rtimes F$, where $F\cong G/P$ is a
cyclic
subgroup of $G$, $|F|\neq 1$ and $(|F|,|P|)=1$. Since $n=|P||F|$, with $P$
and $F$ being cyclic groups of co-prime orders, it follows that 
$\psi(C_n)=\psi(P)\psi(F)$.

Let $Z=C_F(P)$. By Lemma 2.1(6) we have
$$\psi(G)<\psi(P)\psi(F)\biggl(\frac {\psi(Z)}{\psi(F)}+\frac {|P|}{\psi(P)}
\biggr)=
\biggl(\frac {\psi(Z)}{\psi(F)}+\frac {|P|}{\psi(P)}\biggr)\psi(C_n).
\tag{${*}{*}$}$$
Notice that since $P$ is cyclic and $p\geq q+2$, we have
$$\frac  {|P|}{\psi(P)}=\frac {|P|(p+1)}{p|P|^2+1}<\frac {p+1}{p|P|}\leq \frac
{p+1}{p^2}
=\frac 1p+\frac 1{p^2}\leq \frac {q+3}{(q+2)^2}.$$

Consider now the other fraction $\frac {\psi(Z)}{\psi(F)}$. If $C_F(P)=F$, then
$G=P\times F$
is cyclic, a contradiction. So suppose that
$C_F(P)=Z<F$.
Notice that $\psi(F)$ is a product of $\psi(S)$, with $S$ running over all
Sylow
subgroups of $F$. Since $Z$ is a proper subgroup of $F$, also $\psi(Z)$ is a
similar
product, and at least one Sylow subgroup of $Z$, say the Sylow $r$-subgroup of
$Z$, is properly
contained in the Sylow $r$-subgroup of $F$ of order $r^s\geq r$. Hence
$$\frac {\psi(Z)}{\psi(F)}\leq \frac {r^{2(s-1)+1}+1}{r^{2s+1}+1},$$
where $r\geq 3$ and $s\geq 1$.

We claim that the inequality
$$\frac {r^{2(s-1)+1}+1}{r^{2s+1}+1}\leq \frac 1{r^2-r+1}$$
is true.
Indeed, this inequality is equivalent to the following sequence of inequalities:
$$r^{2s+1}+r^2-r^{2s}-r+r^{2(s-1)+1}+1\leq r^{2s+1}+1,$$
$r-rr^{2s-2}-1+r^{2s-2}\leq 0\ $, $r-1\leq (r-1)r^{2s-2}$ and finally $1\leq
r^{2s-2}$, which is true.
Since $r\geq q$, it follows that
$$\frac {\psi(Z)}{\psi(F)}\leq \frac 1{q^2-q+1}.$$

Thus by $({*}{*})$
$$\psi(G)\leq \biggl(\frac  1{q^2-q+1}+\frac {q+3}{(q+2)^2}\biggr)\psi(C_n),$$
and we claim that
$$\frac  1{q^2-q+1}+\frac {q+3}{(q+2)^2}<\frac {((q^2-1)q+1)(q+1)}{q^5+1},$$
which yields a contradiction to our assumptions.
(Notice that this is not true for $q=2$, being equivalent to $341<336$. But
$q\geq 3$,
and for $q=3$ it is barely true, being equivalent to $4087<4375$.)

Multiplying the above inequality by $q^3+1$ and noting that
$q^3+1=(q^2-q+1)(q+1)$, our
claim becomes
$$q+1+\frac {(q+3)(q^3+1)}{q^2+4q+4}<\frac {(q^3-q+1)(q+1)(q^3+1)}{q^5+1}.$$
Now
$$q+1+\frac {(q+3)(q^3+1)}{q^2+4q+4}<q+1+\frac
{(q+3)(q+1)(q^2-q+1)}{q^2+4q+3}=q+1+q^2-q+1=q^2+2,$$
and
$$q^2+2<\frac {(q^3-q+1)(q^4+q^3)}{q^5+1}$$
since this inequality is equivalent to 
$3q^5+q^2+2<q^6+q^3$,
which is true
as $q\geq 3$. So
$$q+1+\frac {(q+3)(q^3+1)}{q^2+4q+4}<\frac {(q^3-q+1)(q^4+q^3)}{q^5+1}<\frac
{(q^3-q+1)(q+1)(q^3+1)}{q^5+1},$$
and our claim follows, yielding a contradiction. The proof of  
Proposition 6 is now complete.
\qed
\enddemo

\heading III. Proof of Proposition 7.
\endheading
In this section we shall prove Proposition 7  which is restated below. 
As previously,
$f(x)$ will denote the function $\frac {((x^2-1)x+1)(x+1)}{x^5+1}$. 
\proclaim {Proposition 7}  Let $q$ be a prime and let $G$ be a  
$q^*$-group of
order $n$. Then
$$\psi(G)= \frac
{((q^2-1)q+1)(q+1)}{q^5+1}\psi(C_n)=f(q)\psi(C_n)\tag{${*}{*}{*}$}$$
if and only if $n=q^2k$ with $(k,q!)=1$ and
$G=(C_q\times C_q)\times C_k$.
\endproclaim
\demo{Proof}
First we show that if $G$ is the group $(C_q\times C_q)\times C_k$
of order $n=q^2k$ with $(k,q!)=1$, then the equality $({*}{*}{*})$
holds. Indeed,
$\psi(G) =\psi(C_q\times C_q)\psi(C_k)=\frac
{(q^2-1)q+1}{\psi(C_{q^2})}\psi(C_n)$ and $\psi(C_{q^2})=\frac {q^5+1}{q+1}$,
so $({*}{*}{*})$ holds, as claimed.

Now we turn to
the "only if" part. So suppose that the $q^*$-group $G$ of order $n$ 
satifies the equality $({*}{*}{*})$. Since $f(q)<1$, $G$ is non-cyclic.
Our aim is to show that $n=q^2k$ with $(k,q!)=1$ and
$$G=(C_q\times C_q)\times C_k.$$

Let $p$ be the maximal prime divisor of $n$.
Our proof is by induction on the size of $p$.

Suppose, first that $p=q$. Then $n=q^r$ and $r\geq 2$. If $r=2$, then $G$ is a
non-cyclic group of order $q^2$. The group $C_q\times C_q$ is the only such
group   and it clearly satisfies both conditions of Proposition 7.
If $q=2$ and $n=2^3$, then it is easy to see that $Q_8$ attains the maximal 
value of $\psi$ among the non-cyclic groups of order $8$. But 
$\psi(Q_8)=27<f(2)\psi(C_8)=\frac 7{11}\cdot 43$, so no non-cyclic group
of order $8$ satisfies $({*}{*}{*})$.
So suppose that either $r=3$ and $q>2$ or $r\geq 4$. Then by 
Theorem 4.4(1) and Lemma 4.2 in [16], the groups $C_q\times C_{q^{r-1}}$ and  $M_{q^r} =  \langle  a, b \ | \ a^{q^r-1}= b^q = 1 , a^b = a^{q^{r-2}+1}\rangle$ attain
 the maximal 
value of $\psi$ among the non-cyclic groups of order $q^r$ and
$$\psi(C_q\times C_{q^{r-1}}) = \psi(M_{q^r})=\frac {q^{2r}+q^3-q^2+1}{q+1}.$$
We claim that 
$$\frac {q^{2r}+q^3-q^2+1}{q+1}<f(q)\psi(C_n)=
\biggl( \frac {q^4+q^3-q^2+1}{q^5+1}\biggr)\cdot \biggl(\frac
{q^{2r+1}+1}{q+1}\biggr)$$
and hence these values for $n$ are impossible. To show this, we need to prove
that
$$(q^{2r}+q^3-q^2+1)(q^5+1)<(q^4+q^3-q^2+1)(q^{2r+1}+1)$$
or
$(q^7+q^4)(q-1)<(q^{2r+3}+q^{2r})(q-1)$ or $q^4<q^{2r}$, which is true, since 
$r\geq 3$. Thus if $n=q^r$ and $G$ satisfies $({*}{*}{*})$, 
then only the case $n=q^2$ and $G=C_q\times C_q$
is possible, as required.

Suppose, now, that $p>q$. By Lemma 2.1(4) $\psi(C_n)\geq \frac q{p+1}n^2$, so
$$\psi(G)= \frac {((q^2-1)q+1)(q+1)}{q^5+1}\psi(C_n)\geq \frac {(q^4+q^3-q^2+1)q}
{(q^5+1)(p+1)}n^2.$$
Hence there exists $x\in G$ such that
$$o(x)>\frac {q^5+q^4-q^3+q}{(q^5+1)(p+1)}n
\quad \text{and}\quad 
[G:\langle x\rangle]< \frac {(q^5+1)(p+1)}{q^5+q^4-q^3+q}.$$

{\bf Case A.} Suppose that $p$ divides $[G:\langle x\rangle]$.

  Then $[G:\langle x\rangle]=pl$ for some positive integer $l$ and
$$pl < \frac {q^5+1}{q^5+q^4-q^3+q}(p+1)<p+1.$$
Hence $l=1$ and $[G:\langle x\rangle]=p$, which implies that
$(q^5+q^4-q^3+q-q^5-1)p<q^5+1$. Hence 
$$p< \frac {q^5+1}{q^4-q^3+q-1},$$
and we claim that 
$$q^5+1<(q^4-q^3+q-1)(q+2)=q^5+q^4-2q^3+q^2+q-2$$
or $2q^3 +3<q^4+q^2+q$, which is true, since $q\geq 2$. Thus $p<q+2$
and since $p>q$, it follows that $p=q+1$. Hence $q=2$, $p=3$ and in particular 
$[G:\langle x\rangle]=3$.

Thus
$$|G|=n=2^a3^b\quad \text{with}\ a,b\geq 1,$$
and $\langle x\rangle$ contains a Sylow $2$-subgroup $Q$ of $G$. Hence
$Q$ is cyclic and the Sylow $3$-subgroup $P$ of $G$ is normal in $G$. Thus
$G=P\rtimes Q$.

Suppose that there exists  $y\in G$ satisfying $[G:\langle y\rangle]=2$. Then 
$\langle y\rangle$ contains $P$, so also $P$ is cyclic. If there exists 
$z\in G\setminus \langle y\rangle$  also satisfying $[G:\langle z\rangle]=2$,
then $\langle y\rangle\cup \{z\}\subseteq C_G(P)$ and $P\leq Z(G)$. But then
$G=P\times Q $ is cyclic, a contradiction. Hence if $z\in G\setminus \langle y\rangle$,
then $o(z)\leq n/3$ and
$\psi(G)\leq \psi(C_{n/2})+(\frac n2)(\frac n3)$. But this inequality implies, 
as shown in the proof of Theorem 1 in [12], that 
$\psi(G)<\frac 7{11}\psi(C_n)=f(2)\psi(C_n)$, in contradiction to our assumptions.

So we may assume that $o(z)\leq n/3$ for all $z\in G$. 
Since $[G:\langle
x\rangle]=3$, it follows that
$$\psi(G)\leq \psi(C_{n/3})+2(\frac n3)(\frac n3)=\psi(C_{2^a})\psi(C_{3^{b-1}})+2(\frac
n3)^2.$$
But this inequality implies, 
as shown in the proof of Theorem 1 in [12], that
$\psi(G)<\frac 7{11}\psi(C_n)=f(2)\psi(C_n)$, in contradiction to our assumptions.

So {\bf Case A} is impossible.

{\bf Case B.} Suppose that $p$ does not divide $[G:\langle x\rangle]$.

Then $\langle x\rangle$ contains a Sylow $p$-subgroup $P$ of $G$. Hence $P$ 
is cyclic and  $\langle x\rangle\leq N_G(P)$. Recall that
$[G:\langle x\rangle]< \frac {(q^5+1)(p+1)}{q^5+q^4-q^3+q}$. Hence
$$[G:\langle x\rangle]=[G:N_G(P)][N_G(P):\langle x\rangle]
=(1+hp)[N_G(P):\langle x\rangle]<\frac {(q^5+1)(p+1)}{q^5+q^4-q^3+q}<p+1$$
for some non-negative integer $h$. It follows that $h=0$ and $P$ is a normal 
cyclic Sylow $p$-subgroup of $G$. Thus, by Lemma 2.1(5), we have
$\psi(G)\leq \psi(P)\psi(G/P)$, with equality if and only if $P$ is central 
in $G$. But
$$\psi(G)= \frac {((q^2-1)q+1)(q+1)}{q^5+1}\psi(C_n)=\frac
{((q^2-1)q+1)(q+1)}{q^5+1}\psi(P)\psi(C_{|G/P|}),$$
so it follows that 
$$\psi(G/P)\geq \frac {((q^2-1)q+1)(q+1)}{q^5+1}\psi(C_{|G/P|}).$$

Suppose, first, that $G/P$ is non-cyclic. Since $G/P$ is a $q^*$-group,
Proposition 6 implies that
$$\psi(G/P)=\frac {((q^2-1)q+1)(q+1)}{q^5+1}\psi(C_{|G/P|})$$ 
and hence 
$\psi(G)=\psi(P)\psi(G/P)$. It follows by Lemma 2.1(5) thet $P\leq Z(G)$
and $G=P\times K$ for some subgroup $K$ of $G$ which is isomorphic to $G/P$.
Since
$$\psi(K)=\frac {((q^2-1)q+1)(q+1)}{q^5+1}\psi(C_{|K|})$$ 
and since 
$q\mid |K|$ and $p\nmid |K|$, it follows by our inductive assumptions that
$K=(C_q\times C_q)\times C_s$, with $(s,q!)=1$. But then 
$G=(C_q\times C_q)\times C_{s|P|}$ with $(s|P|,q!)=1$, as required.

So suppose, finally, that $G/P$ is cyclic. Then $G=P\rtimes F$, where 
$P$ is a cyclic $p$-group and $F$ is a cyclic group of order co-prime to $p$.
Thus $n=|P||F|$ and $\psi(C_n)=\psi(P)\psi(F)$.

Since $G$ is a $q^*$-group satisfying $({*}{*}{*})$, 
Proposition 6 implies that $\psi(G)$ has the second 
largest value among groups of order $n$ and by Lemma 2.1(7) $[F:C_F(P)]$ is a prime.
If $[F:C_F(P)]\neq q$, then a Sylow $q$-subgroup $Q$ of $G$ is contained in $C_F(P)$.
Thus $Q$ is cyclic and it centralizes both $P$ and $F$, so $Q\leq Z(G)$. Hence
$G=Q\times K$, with $|K|=k$, $q\nmid k$ and $K$ is non-cyclic. Now 
$$\psi(G)=\psi(Q)\psi(K)=f(q)\psi(C_n)=f(q)\psi(C_k)\psi(Q)$$
since $Q$ is cyclic. Hence $\psi(K)=f(q)\psi(C_k)$. But since $q\nmid k$, it
follows by Proposition 6 that $\psi(K)\leq f(s)\psi(C_k)$ for some $s>q$. 
As shown in
the Introduction, $f(x)$ is a decreasing function, so $f(s)<f(q)$ and we have 
reached a contradiction.

Hence $[F:C_F(P)]=q$ and $C_F(P)\times P$ is a cyclic subgroup of $G$ of 
index $q$. So there exists $y\in G$ such that
$[G:\langle y\rangle]=q$ and $P\leq \langle y\rangle$.

Suppose, first, that $q=2$. Then $[G:\langle y\rangle]=2$ and 
$P\leq \langle y\rangle$. If there exists $z\in G\setminus \langle y\rangle$ 
such that 
$[G:\langle z\rangle]=2$, then $\langle y\rangle \cup \{z\}\subseteq C_G(P)$.
Hence $P\leq Z(G)$ and $G=P\times F$ is a cyclic group,
a contradiction. So if $z\in G\setminus \langle y\rangle$, then $o(z)\leq n/3$ and
$$\psi(G)\leq \psi(C_{n/2})+(n/2)(n/3).$$ 

If $p=3$, then $n=2^a3^b$ with $a$ and $b$ positive integers, and it
follows by the proof of Theorem 1 in [12] that 
$$\psi(G)<\frac 7{11}\psi(C_n)=f(q)\psi(C_n),$$
in contradiction to our assumptions.

So it remains to deal with two cases: either $q=2$ and $p\geq 5$, or $q\geq 3$. 
Let $Z=C_F(P)$. By Lemma 2.1(6) we have
$$\psi(G)<\psi(P)\psi(F)\biggl(\frac {\psi(Z)}{\psi(F)}+\frac {|P|}{\psi(P)}\biggr)=
\biggl(\frac {\psi(Z)}{\psi(F)}+\frac {|P|}{\psi(P)}\biggr)\psi(C_n).$$
Since $P$ is cyclic, we have
$$\frac  {|P|}{\psi(P)}=\frac {|P|(p+1)}{p|P|^2+1}<\frac {p+1}{p|P|}\leq \frac
{p+1}{p^2}
=\frac 1p+\frac 1{p^2}.$$
Consider now the other fraction $\frac {\psi(Z)}{\psi(F)}$. As shown in the 
proof of Proposition 6, $[F:Z]=q$ implies that
$$\frac {\psi(Z)}{\psi(F)}\leq \frac 1{q^2-q+1}.$$

If $q=2$ and $p\geq 5$, then $p\geq q+3$. Thus
$$\frac  {|P|}{\psi(P)}<\frac 1p+\frac 1{p^2}\leq \frac {q+4}{(q+3)^2}
=\frac 6{25}\quad \text{and}\quad 
\frac {\psi(Z)}{\psi(F)}\leq \frac 1{q^2-q+1}=\frac 13.$$
Hence
$$\psi(G)<(\frac 13+\frac 6{25})\psi(C_n)=\frac {43}{75}\psi(C_n)
<\frac 7{11}\psi(C_n)=f(q)\psi(C_n),$$
which contradicts our assumptions.

If $q>2$, then $q\geq 3$ and $p\geq q+2$. Thus
$$\frac  {|P|}{\psi(P)}<\frac 1p+\frac 1{p^2}\leq \frac {q+3}{(q+2)^2}
\quad \text{and}\quad 
\frac {\psi(Z)}{\psi(F)}\leq \frac 1{q^2-q+1}.$$
Hence
$$\psi(G)<\biggl(\frac 1{q^2-q+1}+\frac {q+3}{(q+2)^2}\biggr)\psi(C_n)$$
and as shown in the proof of Proposition 6, this inequality implies that
$$\psi(G)<\frac {
((q^2-1)q+1)(q+1)}{q^5+1}\psi(C_n)=f(q)\psi(C_n),$$
which contradicts our assumptions. 
 
So in {\bf Case B} $P$ is normal in $G$, and only groups with $G/P$ 
non-cyclic can satisfy condition  $({*}{*}{*})$. As shown above, these
groups are as required.

The proof of Proposition 7 is now complete.

\qed
\enddemo

\bigskip

\Refs

\ref
\no 1
\by H. Amiri, S.M. Jafarian Amiri, I.M. Isaacs
\paper Sums of element orders in finite groups
\jour Comm. Algebra 
\vol 37 
\yr 2009
\pages 2978-2980
\endref
\ref
\no 2
\by H. Amiri, S.M. Jafarian Amiri
\paper Sums of element orders on finite groups of the same order
\jour J. Algebra Appl.
\vol 10 (2)
\yr 2011
\pages187-190
\endref
\ref
\no 3
\by H. Amiri, S.M. Jafarian Amiri
\paper Sum of element orders of maximal subgroups of the symmetric group
\jour Comm. Algebra
\vol 40 (2)
\yr 2012
\pages 770-778
\endref
\ref
\no 4
\by M. Baniasad Asad, B. Khosravi
\paper A Criterion for Solvability of a Finite Group by the Sum of Element Orders
\jour J. Algebra
\vol 516
\yr 2018
\pages 115-124
\endref 
\ref
\no 5
\by S.M. Jafarian Amiri
\paper Second maximum sum of element orders on finite nilpotent groups
\jour Comm. Algebra
\vol 41 (6)
\yr 2013
\pages 2055-2059
\endref
\ref
\no 6
\by S.M. Jafarian Amiri
\paper Maximum sum of element orders of all proper subgroups of $PGL(2,q)$
\jour Bull. Iran. Math. Soc.
\vol 39 (3) 
\yr 2013
\pages 501-505
\endref
\ref
\no 7
\by S.M. Jafarian Amiri
\paper Characterization of $A_5$ and $PSL(2,7)$ by sum of element orders
\jour Int. J. Group Theory
\vol 2 (2)
\yr 2013
\pages 35-39
\endref
\ref
\no 8
\by S.M. Jafarian Amiri, M. Amiri
\paper Second maximum sum of element orders on finite  groups
\jour J. Pure Appl. Algebra
\vol 218 (3)
\yr 2014
\pages 531-539
\endref
\ref
\no 9
\by S.M. Jafarian Amiri, M. Amiri
\paper Sum of the products of the orders of two distinct elements in finite groups
\jour Comm. Algebra
\vol 42 (12)
\yr 2014
\pages  5319-5328
\endref
\ref
\no 10
\by S.M. Jafarian Amiri, M. Amiri
\paper Characterization of $p$-groups by sum of the element orders
\jour Publ. Math. Debrecen
\vol 86 (1-2)
\yr 2015
\pages 31-37
\endref
\ref
\no 11
\by S.M. Jafarian Amiri, M. Amiri
\paper Sum of the Element Orders in Groups with the Square-Free Order 
\jour Bull. Malays. Math. Sci. Soc.
\vol 40
\yr 2017
\pages  1025-1034
\endref
\ref
\no 12
\by M. Herzog, P. Longobardi, M. Maj
\paper An exact upper bound for sums of element orders in non-cyclic finite groups
\jour J. Pure Appl. Algebra
\vol 222 (7) 
\yr 2018
\pages 1628-1642
\endref
\ref
\no 13
\by M. Herzog, P. Longobardi, M. Maj
\paper  Two new criteria for solvability of finite groups in finite groups
\jour J. Algebra
\vol 511
\yr 2018
\pages 215-226
\endref
\ref
\no 14
\by M. Herzog, P. Longobardi and M. Maj
\paper Sums of element orders in groups of order $2m$ with $m$ odd
\jour Comm. Algebra
\toappear
\endref
\ref
\no 15
\by Y. Marefat, A. Iranmanesh, A. Tehranian
\paper On the sum of element orders of finite simple groups
\jour J. Algebra Appl.
\vol12 (7)
\yr 2013
\pages 135-138
\endref
\ref
\no 16
\by R. Shen, G. Chen and C. Wu 
\paper On groups with the second largest value of the sum of element orders
\jour Comm. Algebra
\vol 43 (6)
\yr 2015
\pages 2618-2631
\endref
\ref
\no 17
\by M. T\u{a}rn\u{a}uceanu, D.G. Fodor
\paper On the sum of element orders of finite abelian groups
\jour Sci. An. Univ. "A1.I. Cuza" Iasi, Ser. Math.
\vol LX
\yr 2014
\pages 1-7
\endref

\endRefs
\end{document}